# STATIONARY DISTRIBUTIONS OF MULTI-TYPE TOTALLY ASYMMETRIC EXCLUSION PROCESSES

By Pablo A. Ferrari and James B. Martin

*Universidade de São Paulo and University of Oxford*

We consider totally asymmetric simple exclusion processes with $n$ types of particle and holes ($n$-TASEPs) on $\mathbb{Z}$ and on the cycle $\mathbb{Z}_N$. Angel recently gave an elegant construction of the stationary measures for the 2-TASEP, based on a pair of independent product measures. We show that Angel's construction can be interpreted in terms of the operation of a discrete-time $M/M/1$ queueing server; the two product measures correspond to the arrival and service processes of the queue. We extend this construction to represent the stationary measures of an $n$-TASEP in terms of a system of queues in tandem. The proof of stationarity involves a system of $n$ 1-TASEPs, whose evolutions are coupled but whose distributions at any fixed time are independent. Using the queueing representation, we give quantitative results for stationary probabilities of states of the $n$-TASEP on $\mathbb{Z}_N$, and simple proofs of various independence and regeneration properties for systems on $\mathbb{Z}$.

**1. Introduction.** Consider first the *totally asymmetric simple exclusion process* (TASEP) on $\mathbb{Z}$ or on the cycle $\mathbb{Z}_N$. Some sites $i \in \mathbb{Z}$ (or $\mathbb{Z}_N$) contain a particle; at the others we say there is a hole. The dynamics of the system are as follows: a bell rings as a Poisson process of rate 1 at each site independently; when the bell at site $i$ rings, if there is a particle at site $i$ and a hole at site $i-1$, they exchange. Put another way, each particle in the system tries to jump to the left at rate 1; the jump succeeds whenever the site to its left is unoccupied. See, for example, [9, 10] for a wealth of information on the TASEP and its relatives.

We will denote a configuration of the TASEP by $u = (u(j), j \in \mathbb{Z} \text{ or } \mathbb{Z}_N) \in \mathcal{U}_1$ or $\mathcal{U}_1^{(N)}$, where $\mathcal{U}_1 = \{1, \infty\}^{\mathbb{Z}}$ and $\mathcal{U}_1^{(N)} = \{1, \infty\}^{\mathbb{Z}_N}$. We set $u(i) = 1$ if there is a particle at $i$ and $u(i) = \infty$ if there is a hole at $i$. (This notation is not standard but will be convenient later.)

The stationary measures for the TASEP are as follows. On $\mathbb{Z}_N$, the number of particles in the system, $p$, is conserved by the dynamics. For fixed $p$, the dynamics are given by an irreducible Markov chain with finite state-space, whose stationary measure is uniform over all $\binom{N}{p}$ possible states. Any stationary measure for the system is some linear combination of these uniform distributions.

On $\mathbb{Z}$ there is a one-parameter family of translation-invariant extremal stationary measures $\nu_\lambda$, $\lambda \in [0, 1]$. Under $\nu_\lambda$, each site $i$ is occupied by a particle independently and with probability $\lambda$; the measure $\nu_\lambda$ is the *product measure* with density $\lambda$ of particles (or the measure of a *Bernoulli process with rate $\lambda$*). The only other extremal stationary measures are the so-called *blocking measures*: the measure concentrated on the configuration with particles at all negative sites and holes at all nonnegative sites, and its translates.

We now consider a *TASEP with $n$ classes of particle* or *$n$-TASEP*, for $n \geq 1$. Now a configuration $u = (u(j), j \in \mathbb{Z}$ or $\mathbb{Z}_N)$ of the system is a member of $\mathcal{U}_n$ or $\mathcal{U}_n^{(N)}$, where

$$\mathcal{U}_n = (\{1, 2, \ldots, n\} \cup \{\infty\})^{\mathbb{Z}}, \qquad \mathcal{U}_n^{(N)} = (\{1, 2, \ldots, n\} \cup \{\infty\})^{\mathbb{Z}_N}.$$

If $u(i) = r \leq n$ we say that there is a particle of class $r$ at site $i$, while if $u(i) = \infty$ we say that there is a hole at $i$.

The dynamics are as follows. Bells ring at rate 1 at each site independently. When the bell at site $i$ rings, the values at $i-1$ and $i$ swap if $u(i-1) > u(i)$, and remain unchanged otherwise. That is, the process jumps from $u$ to the configuration $u^{(i-1,i)}$ where

$$u^{(i-1,i)}(j) = u(j) \qquad \text{if } j \notin \{i-1, i\},$$

$$u^{(i-1,i)}(i-1) = \min\{u(i-1), u(i)\},$$

$$u^{(i-1,i)}(i) = \max\{u(i-1), u(i)\}.$$

Put another way, each particle in the system tries to jump to the left at rate 1. The jump succeeds when the site to its left contains a hole or contains a particle with a higher class (in which case the two particles exchange positions).

(Note that the way we have defined it, the $n$-type TASEP has $n + 1$ possible values for a site: $n$ types of particle and holes. It might be more elegant simply to regard holes as particles of type $n + 1$, and to replace the value $\infty$ by $n + 1$ throughout; but the distinguished state $\infty$ will be convenient later for various reasons.)

One natural way in which $n$-type TASEPs arise is from couplings of several 1-type TASEPs. Let $\eta_t^{(1)}$, $\eta_t^{(2)}, \ldots, \eta_t^{(n)}$, $t \geq 0$ be $n$ processes realizing a TASEP, started from initial conditions such that $\eta_0^{(1)}(j) \geq \eta_0^{(2)}(j) \geq \cdots \geq \eta_0^{(n)}(j)$ for all $j$. (I.e., whenever there is a particle at site $j$ in process $m$, there is also one in the processes $m+1, m+2, \ldots, n$.) Suppose one couples the processes by using the same processes of bells at each site for all of them. Then the ordering $\eta_t^{(1)}(j) \geq \eta_t^{(2)}(j) \geq \cdots \geq \eta_t^{(n)}(j)$ continues to hold for all $t$ (this is an instance of the *basic coupling* of Liggett [8, 9]). Let $u_t(j) = \inf\{m : \eta_t^{(m)}(j) = 1\}$ [with $u_t(j) = \infty$ if $\eta_t^{(m)}(j) = \infty$ for all $m$]. Then $u_t$ realizes an $n$-type TASEP.

The 2-TASEP has been studied from several perspectives. The existence of a translation-invariant stationary measure for the process on $\mathbb{Z}$, with densities $\lambda_1$ and $\lambda_2$ of first- and second-class particles (for $0 < \lambda_1 < \lambda_1 + \lambda_2 < 1$), was proved by Liggett [8] to demonstrate ergodic properties of the TASEP. The uniqueness and extremality of this measure were shown by Ferrari, Kipnis and Saada [5] and Speer [12] (the only other extremal invariant distributions are those concentrated on a single state—"blocking measures" as for the 1-TASEP above). Derrida, Janowsky, Lebowitz and Speer [2] and Speer [12] construct the measure explicitly (both in finite and infinite volume) using a matrix method and show various regeneration and asymptotic properties; Ferrari, Fontes and Kohayakawa [4] give probabilistic interpretations and proofs of the measure and its properties.

Recently Angel [1] gave an elegant construction of this stationary measure based on two independent product measures with densities $\lambda_1$ and $\lambda_1 + \lambda_2$ (and an analogous construction for the case of $\mathbb{Z}_N$); the proof involves providing bijections between certain families of binary trees and pairs of binary sequences.

We will show that Angel's construction can be rewritten in terms of the operation of a queueing server (namely, an $M/M/1$ queue in discrete time). The two independent product measures correspond to the arrival process and the service process of the queue. The stationary measure for the 2-TASEP corresponds to the output process of the queue. Sites of the TASEP are interpreted as times in the queueing process. If the queue has a departure at time $i$, then one puts a first class particle at site $i$. If the queue has an "unused service" at time $i$ (i.e., a service is available but there is no customer present) then one puts a second-class particle at site $i$. If there is no service available at time $i$, then there is a hole at site $i$. (Full descriptions are given in the next section.)

We then generalize this result to the $n$-TASEP with $n \geq 3$. First a remark about the set of stationary measures. Let $\lambda_1, \ldots, \lambda_n \in (0, 1)$ with $\sum_{r=1}^n \lambda_r < 1$. For the process on $\mathbb{Z}$, there is a unique translation-invariant stationary measure with density $\lambda_r$ of particles of type $r$, for each $r$. These stationary measures are extremal, and the only other stationary measures are blocking measures, concentrated on a single configuration. These facts can be proved almost identically to the proofs given in, for example, [5] and [12] for the case $n = 2$. The case of the $n$-TASEP on $\mathbb{Z}_N$ is analogous; for any $p_1, \ldots, p_n \in \{0, 1, \ldots, N\}$ with $\sum_{r=1}^{n} p_n \leq N$, there is an extremal stationary measure concentrated on configurations with exactly $p_r$ particles of type $r$, for each $r$; these are the only extremal stationary measures.

We construct a representation of these stationary measures based on a system of $n - 1$ *queues in tandem*. Take the case of $\mathbb{Z}$, say. Consider a set of $n$ independent Bernoulli processes (product measures), such that the $m$th has density $\lambda_1 + \cdots + \lambda_m$, for each $m$. The service process of the $m$th queue is given by the $(m + 1)$st product measure. The first product measure acts as the arrival process to the first queue. Thereafter, the arrival process to the $m$th queue, $m > 1$, is given by the output process from the $(m - 1)$st queue. The $m$th queue is a "priority queue" with $m$ types of customer; at a service time, the customer who departs is the one whose class number is lowest out of those present. We will show that the distribution of the output process of the $(n - 1)$st queue provides the stationary distribution of the $n$-TASEP with densities $\lambda_1, \ldots, \lambda_n$. (We describe the operation of the system fully in the next section.)

The proof of stationarity has the following structure. First we define dynamics on a set of $n$ configurations of particles on the line (the "multiline process"). On each line, the local transitions are those of a TASEP, and the lines are coupled in such a way that the bells on neighboring lines always ring at the same time, either at the same site or at neighboring sites. We will show that a collection of $n$ independent product measures is stationary under these dynamics. (In fact, it will turn out that in this equilibrium, the multiline process provides a coupling of $n$ TASEP processes, whose marginal distributions at any fixed time are independent.) The second step of the proof is to show that, under these dynamics, the final output process derived as in the previous paragraph is realizing an $n$-TASEP.

In the case $n = 2$, the two lines of the multiline process correspond to the arrival and departure processes of the queue. (We note that, especially in the case $n = 2$, the structure of the proof is similar to that used in recent work by Duchi and Schaeffer [3], where a related two-line process is used to study combinatorial aspects of TASEPs with one or two types of particle, on $\mathbb{Z}_N$ or on a finite strip with boundaries.)

In Section 2, we give precise definitions and state the main results, first for $n = 2$ (Angel's construction and the queueing interpretation) and then for general $n$. In Section 3, the first step of the proof is carried out: the "multiline process" is defined and it is proved that a collection of independent product measures is stationary under its dynamics. (For the case of $\mathbb{Z}_N$, the analogous stationary distribution is a collection of independent uniform distributions over configurations with a fixed number of particles.) The second step of the proof described in the previous paragraph is done in Section 4.

In Section 5 we use the multiline representation to prove a quantitative result about the probabilities of states in the stationary distribution of an $n$-TASEP on the ring $\mathbb{Z}_N$. Let $p_1, \ldots, p_n \in \{0, 1, \ldots, N\}$ with $\sum_{r=1}^{n} p_r \leq N$ and consider the stationary distribution concentrated on configurations with exactly $p_r$ particles of class $r$ for each $r$. Then the probability of every state is an integer multiple of $M^{-1}$, where

$$(1.1) \qquad M = \left[ \binom{N}{p_1} \binom{N}{p_1 + p_2} \cdots \binom{N}{p_1 + p_2 + \cdots + p_n} \right].$$

This confirms a conjecture of Angel [1]. (See also the discussion by Mallick, Mallick and Rajewsky in [11], where the matrix method is used to investigate the case $n = 3$.) We can then identify those states with probability $M^{-1}$ (the smallest possible). A configuration $u$ is minimal in this way if the following holds for each $j \in \mathbb{Z}_N$: if $u(j) = \infty$ then $u(j + 1) = \infty$ or $u(j + 1) = n$, while if $u(j) = m < \infty$ then $u(j + 1) \geq m - 1$.

For example, consider $n = 2$, $N = 9$, $p_1 = 3$, $p_2 = 3$. Possible "minimal states" include: $(\infty, \infty, \infty, 2, 2, 2, 1, 1, 1)$ and $(\infty, 2, \infty, \infty, 2, 1, 2, 1, 1)$.

For the extreme case of $N$ types of particle on $\mathbb{Z}_N$, with $p_m = 1$ for all $m = 1, \ldots, N$, the only minimal states are $(N, N-1, \ldots, 3, 2, 1)$ and its cyclic shifts (giving $N$ minimal states in all). In this case one has

$$M = \frac{(N!)^{N-1}}{(2! \, 3! \cdots (N-1)!)^2}.$$

An important motivation for much of the previous work on the 2-TASEP on $\mathbb{Z}$ was its application to the study of shock measures for the 1-TASEP; see, for example, [2, 4, 5, 12]. In Section 6, we show how several results developed in this context—in particular, concerning the "stationary measure as seen from a second-class particle"—can be proved very simply, and often strengthened, using the queueing representation. We also give extensions of some of these results to the case of the $n$-TASEP for $n > 2$.

In passing we mention another paper [6] in which we discuss (partly in survey form) results related to those in this paper, but with a somewhat different approach. For example, a proof that product measure is invariant for the multiclass TASEP on $\mathbb{Z}$ is described much more directly in terms of reversibility arguments in infinite volume, in contrast to the combinatorial approach in terms of bijections on the state-space of the process on $\mathbb{Z}_N$, used in Section 3 of this paper. The aim of [6] is partly to illustrate the key role played by reversibility properties, and also to show that the same multiclass measures are stationary for a variety of underlying dynamics (including, e.g., a multiclass discrete-space version of the Aldous–Diaconis–Hammersley process and certain discrete-time versions of the multiclass TASEP); the method of proof is described independently of the dynamics, as far as possible. The proof in Sections 3 and 4 of this paper is much more specific to the TASEP, but (perhaps unlike the approach of [6]) adapts quite naturally to the case of the process on $\mathbb{Z}_N$.

## 2. Main result.

2.1. $n = 2$: *Angel's construction and the queueing interpretation.* We first describe Angel's "collapse" process for the case of $\mathbb{Z}_N$.

Choose $p_1$ and $p_2$ with $0 < p_1 < p_1 + p_2 < N$. Let $\mathcal{A}$ be a subset of $\{1, 2, \ldots, N\}$ of size $p_1$, chosen uniformly at random, and let $\mathcal{S}$ be another subset, of size $p_1 + p_2$, chosen uniformly at random and independently of $\mathcal{A}$.

We will construct a set $\mathcal{D} \subseteq \mathcal{S}$ as follows. $\mathcal{D}$ is initially empty. Consider each site $i \in \mathcal{A}$ in turn (in an arbitrary order). Now we add to the set $\mathcal{D}$ the first element in the list $\{i, i+1, i+2, \ldots\}$ which is in $\mathcal{S}$ and has not yet been added to $\mathcal{D}$. (We "move $i$ to the right," until we find an available space for it; the set $\mathcal{S}$ is the set of spaces which are available to be used.) Notice that a space $j \in \mathcal{S}$ is "filled"—that is, contained in $\mathcal{D}$—if and only if there is some interval $[i, j]$ which contains at least as many members of $\mathcal{A}$ as it does members of $\mathcal{S}$ (since then, by moving the members of $\mathcal{A}$ to the right until they encounter a free space, one of them will inevitably be allocated to the space $j$); so indeed, the order in which the sites of $\mathcal{A}$ are considered does not affect the set $\mathcal{D}$ which is finally obtained by this procedure. Thus

$$(2.1) \qquad j \in \mathcal{D} \qquad \text{iff } j \in \mathcal{S} \text{ and, for some } i, |\mathcal{A} \cap [i, j]| \geq |\mathcal{S} \cap [i, j]|.$$

The set $\mathcal{D}$ constructed in this way (of size $p_1$) gives the locations of first-class particles. The set $\mathcal{S} \setminus \mathcal{D}$ (of size $p_2$) gives the locations of second-class particles. The set $\mathcal{S}^c$ (of size $N - p_1 - p_2$) gives the locations of holes.

The distribution of the configuration generated in this way is the stationary distribution of the TASEP with 2 types of particle on $\mathbb{Z}_N$, with $p_1$ and $p_2$ particles of first and second class, respectively, [1].

In later notation, we will write, for $i \in \mathbb{Z}_N$, $a(i) = 1$ if $i \in \mathcal{A}$ and $\infty$ otherwise, $s(i) = 1$ if $i \in \mathcal{S}$ and $\infty$ otherwise, and $d(i) = 1$ if $i \in \mathcal{D}$ and $\infty$ otherwise. Thus $a$, $s$ and $d$ are configurations from the set $\mathcal{U}_1^{(N)} = \{1, \infty\}^{\mathbb{Z}_N}$ defined earlier. If $d$ is derived from $a$ and $s$ by the process described above we will write $d = F_1^{(N)}(a, s)$.

The case of $\mathbb{Z}$ is analogous. Now we choose $\lambda_1$ and $\lambda_2$ with $0 < \lambda_1 < \lambda_1 + \lambda_2 < 1$. For each $i \in \mathbb{Z}$, let $i \in \mathcal{A}$ with probability $\lambda_1$, independently for each $i$. For each $i \in \mathbb{Z}$, let $i \in \mathcal{S}$ with probability $\lambda_1 + \lambda_2$, independently for each $i$ and independently of the set $\mathcal{A}$. As in (2.1) we can then define

$$(2.2) \quad j \in \mathcal{D} \qquad \text{iff } j \in \mathcal{S} \text{ and, for some } i \leq j, |\mathcal{A} \cap [i, j]| \geq |\mathcal{S} \cap [i, j]|.$$

Equivalently (but rather less formally now that $\mathcal{A}$ is infinite), $\mathcal{D}$ can be constructed by taking each site in $\mathcal{A}$ in turn (in an arbitrary order), and moving it to the right until we find an available space in $\mathcal{S}$ for it.

Again, the configuration of first-class particles (sites of $\mathcal{D}$), second-class particles (sites of $\mathcal{S} \setminus \mathcal{D}$), and holes (sites of $\mathcal{S}^c$) is a sample from the stationary distribution of the TASEP with 2 types of particle on $\mathbb{Z}$, with densities $\lambda_1$ and $\lambda_2$ of first- and second-class particles, respectively [1].

As above, we will write $a(i) = 1$ if $i \in \mathcal{A}$ and $\infty$ otherwise, $s(i) = 1$ if $i \in \mathcal{S}$ and $\infty$ otherwise, and $d(i) = 1$ if $i \in \mathcal{D}$ and $\infty$ otherwise, where now $i \in \mathbb{Z}$. Thus $a$, $s$ and $d$ are configurations from $\mathcal{U}_1 = \{1, \infty\}^{\mathbb{Z}}$. Note that by construction the random variables $a(i), i \in \mathbb{Z}$ and $s(i), i \in \mathbb{Z}$ are all independent. If $d$ is derived from $a$ and $s$ by the process described above we will write $d = F_1(a, s)$.

We now explain the interpretation in terms of a queueing server for the case of $\mathbb{Z}$.

The sites $i \in \mathbb{Z}$ represent times. The set $\mathcal{A}$ represents the set of arrival times, when a (single) customer arrives at the queue. The set $\mathcal{S}$ is the set of service times of the queue. At such times it is possible for a (single) customer to depart. The queue thus has *arrival rate* $\lambda_1$ and *service rate* $\mu = \lambda_1 + \lambda_2$. At each time in $\mathcal{S}$, a customer departs if any customer is present in the queue (including one that has just arrived at that same moment). Then $\mathcal{D}$ is the set of times at which a customer departs from the queue. $\mathcal{S} \setminus \mathcal{D}$ is the set of *unused service times*.

We define

$$(2.3) \qquad Q_j = \left[ \sup_{i \le j} \{ |\mathcal{A} \cap [i, j]| - |\mathcal{S} \cap [i, j]| \} \right]_+.$$

Note that since the service rate is larger than the arrival rate, the law of large numbers yields that the right-hand side of (2.3) is finite with probability 1. From the definition (2.3), and recalling that $a(i) = 1$ [resp. $s(i) = 1$] iff $i \in \mathcal{A}$ (resp. $i \in \mathcal{S}$), we can deduce a set of "recurrences," namely that

$$(2.4) \qquad Q_j = [Q_{j-1} + I(a(j) = 1) - I(s(j) = 1)]_+$$

for each $j \in \mathbb{Z}$. Thus $Q_j$ may be interpreted as the *queue-length* after time $j$ (since (2.4) has the following interpretation: at each time-step, the queue-length increases by 1 if there is an arrival and no service, decreases by 1 if it is nonzero and there is a service and no arrival, and otherwise remains the same). An equivalent formulation of (2.2) above is that

$$(2.5) \qquad d(j) = 1 \qquad \text{iff } s(j) = 1 \text{ and either } Q_{j-1} > 0 \text{ or } a(j) = 1.$$

From (2.3), $Q_j$ is a function of the variables $(a(i), s(i), i \le j)$. Since the variables $a(i), s(i), i \in \mathbb{Z}$ are all independent, (2.4) shows that, conditional on the value of $Q_j$, $Q_{j+1}$ is independent of $(Q_r, r < j)$, so that $(Q_j, j \in \mathbb{Z})$ is a discrete-time Markov chain. In fact it is a birth-and-death chain with one-step transition probabilities given by $i \to i + 1$ with probability $\lambda_1(1 - \lambda_1 - \lambda_2)$, $i \to i - 1$ with probability $I(i > 0)(1 - \lambda_1)(\lambda_1 + \lambda_2)$ and otherwise $i \to i$. The stationary distribution of this chain is geometric with parameter $\lambda_1/(\lambda_1 + \lambda_2)$. From (2.3) and the stationarity of $a$ and $s$, the process $(Q_j, j \in \mathbb{Z})$ is also stationary [and so $Q_j \sim \text{Geom}(\lambda_1/(\lambda_1 + \lambda_2))$ for all $j$].

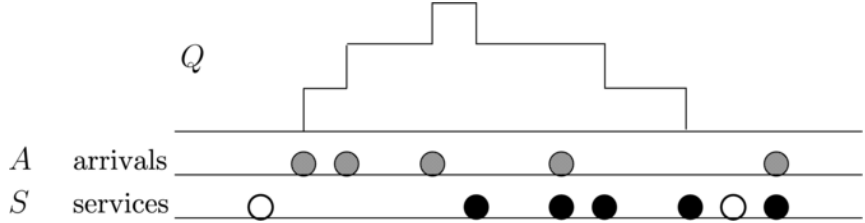

Fig. 1. *Queueing representation of the two-type measure. Arrivals are represented by grey particles. Used services (first-class particles) are represented by black particles and unused ones (second-class particles) by white particles.*

A form of *Burke's theorem* for the discrete-time $M/M/1$ queue [7] applies here. Since the arrival process $\mathcal{A}$ is a Bernoulli process with rate $\lambda_1$, and the service process $\mathcal{S}$ is a Bernoulli process (independent of the arrivals) with rate $\mu > \lambda_1$, the departure process $\mathcal{D}$ is itself a Bernoulli process of rate $\lambda_1$.

The queueing formalism can also be used to describe the process on $\mathbb{Z}_N$. Define

$$(2.6) \qquad Q_j = \left[ \sup_{i \in \mathbb{Z}_N} \{ |\mathcal{A} \cap [i, j]| - |\mathcal{S} \cap [i, j]| \} \right]_+.$$

Because $|\mathcal{A}| < |\mathcal{S}|$, the sup on the RHS of (2.6) is attained by some $i$ which satisfies $|[i, j]| < N$; that is, $i \neq j + 1$ and so the interval $[i, j]$ is not the whole of $\mathbb{Z}_N$. Then, just as before, one can again deduce the relations in (2.4) for each $j \in \mathbb{Z}_N$, and the equivalence of (2.1) and (2.5) also follows. In terms of the collapse process one could give the following description: each particle starts at some location $i$ in the process $\mathcal{A}$ and is allocated some location $k$ in the process $\mathcal{S}$; then $Q_j$ is the number of particles for which the corresponding $i$ and $k$ satisfy $i \neq k$ and $j \in [i, k - 1]$.

2.2. *Extension to several classes of particle/customer.* We now explain how to extend this framework to construct the stationary distribution of the $n$-TASEP for any $n$. In the queueing context, this will correspond to the output process of a series of $n - 1$ queues in tandem with various classes of customer; the output process of one queue acts as the input to the next.

We begin with the case of $\mathbb{Z}$. We first define the behavior of a queueing server whose arrival process contains $m$ classes of arrival (and holes) and whose departure process contains $m + 1$ classes of departure (along with holes).

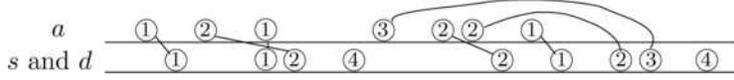

FIG. 2. *Multi-type representation for 3 types of arrivals. Each arrival has been linked to the corresponding service. Note: in this figure and in Figure 3 we have used a "FIFO" (first-in first-out) discipline between customer of the same class; when two or more customers present have joint highest priority, the one that departs is the one that arrived earliest. However, this choice of queueing discipline is entirely arbitrary since for the purposes of the TASEP construction, all customers of the same class may be presumed identical. For example, the second and the third of the 2nd-class customers in the diagram could depart in the opposite order; the departure process d is unchanged.*

Let $a = (a(i), i \in \mathbb{Z}) \in \mathcal{U}_m$ be the arrival process, with $a(i) \in \{1, 2, \ldots, m\} \cup \{\infty\}$. The value $a(i) = r < \infty$ indicates that a customer of class $r$ arrives at time $i$, while $a(i) = \infty$ means that no customer arrives at time $i$.

Let $s = (s(i), i \in \mathbb{Z}) \in \mathcal{U}_1$ be the service process, with $s(i) \in \{1, \infty\}$. A service is available at time $i$ if $s(i) = 1$, and no service is possible if $s(i) = \infty$.

We will derive a departure process $d = (d(i), i \in \mathbb{Z}) \in \mathcal{U}_{m+1}$ from $a$ and $s$, with $d(i) \in \{1, 2, \ldots, m+1\} \cup \{\infty\}$.

Customers of class 1 have the highest priority, and a customer of class $r$ has higher priority than one of class $s$ when $r < s$.

If there is a service possible at time $i$ [i.e., $s(i) = 1$], then the highest priority customer of those currently in the queue, including any customer who has just arrived at time $i$ itself, departs. (If there is currently more than one customer of jointly highest priority, one of these departs.) If the departing customer has class $k$, then $d(i) = k$. If there is no customer in the queue [in particular this implies $a(i) = \infty$], then an "unused service" occurs, and we put $d(i) = m + 1$. (One could imagine that there is a limitless supply of customers of class $m + 1$ available, who get served only when the server would otherwise be idle.) If there is no service at time $i$ [i.e., $s(i) = \infty$] then put $d(i) = \infty$. No departure occurs; all customers in the queue, along with any customer arriving at $i$, remain in it. Note that $d(i) = \infty$ iff $s(i) = \infty$.

In the situations we need, the arrival and departure processes will be stationary and independent, and the long-run intensity of arrivals to the queue will be lower than that of services. In this case, a stationary evolution for the queue, corresponding to the verbal description just given, can be defined in a standard way; see Section 2.3 below.

We define a function (or "queueing operator") $F_m : \mathcal{U}_m \times \mathcal{U}_1 \mapsto \mathcal{U}_{m+1}$, representing the operation of this queue, by $F_m(a, s) = d$.

Now fix $n > 1$. Let $\mathcal{X}$ be the state-space, $\{1, \infty\}^{\mathbb{Z} \times \{1, \ldots, n\}} = \mathcal{U}_1^{\{1, 2, \ldots, n\}}$ representing configurations of particles on $n$ lines. For $x = (x_m(i), i \in \mathbb{Z}, m \in \{1, 2, \ldots, n\}) \in \mathcal{X}$, we say that $x$ has a particle at site $i$ on line $m$ if $x_m(i) = 1$, and that $x$ has a hole at site $i$ on line $m$ if $x_m(i) = \infty$.

Let $\lambda_r \in (0,1)$ for $r \in \{1, 2, \ldots, n\}$, with $\lambda_1 + \cdots + \lambda_n < 1$. We consider the product distribution on the space $\mathcal{X}$ under which all the $x_m(i)$ are independent, with $\mathbb{P}(x_m(i) = 1) = \lambda_1 + \cdots + \lambda_m$.

Consider a richer state-space $\mathcal{V}$ consisting of configurations $v = (v_m(i)$, $i \in \mathbb{Z}$, $m \in \{1, 2, \ldots, n\})$ with $v_m(i) \in \{1, 2, \ldots, m\} \cup \{\infty\}$. If $v_m(i) = \infty$ we say that $v$ has a hole at site $i$ on line $m$, while if $v_m(i) = r < \infty$ then $v$ has a particle of class $r$ at site $i$ on line $m$. Thus $v_m$, the $m$th line of $v$, is a configuration in $\mathcal{U}_m$; All the particles on line $m$ have classes in $\{1, 2, \ldots, m\}$.

Given a configuration $x \in \mathcal{X}$, we will associate a configuration $v = Vx \in \mathcal{V}$ to it. We will have $v_m(i) = \infty$ iff $x_m(i) = \infty$; that is, the particles in $v$ and $x$ are in exactly the same places, and deriving $v$ from $x$ corresponds to assigning classes to the particles of $x$.

First set $v_1(\cdot) = x_1(\cdot)$. (All the particles on the top line have class 1.) Now recursively, define $v_{m+1}(\cdot) = F_m(v_m(\cdot), x_{m+1}(\cdot))$ for $m = 1, 2, \ldots, n-1$, where $F_m$ is the multiclass queueing operator defined above. Hence we have a sequence of $n-1$ queues: the $m$th queue has arrival process $v_m$ and service process $x_{m+1}$, and outputs the departure process $v_{m+1}$ (which is then used as the arrival process for queue $m + 1$).

Note that the intensity of services on line $m$ is $\lambda_1 + \cdots + \lambda_m$, which is increasing in $m$. Applying Burke's theorem repeatedly, one obtains that for $r \leq m \leq n$, the set of particles of class $r$ or lower on line $m$ is a Bernoulli process, and in particular has rate (or density) $\lambda_1 + \cdots + \lambda_r$. The set of particles of class *exactly* $r$ on line $m$ thus has density $\lambda_r$ (in the sense that each site has probability $\lambda_r$ of containing such a particle), but may not be a Bernoulli process.

THEOREM 2.1. *The distribution of $v_m$ (the $m$th line of the configuration $v$) is the stationary distribution of the $m$-type TASEP on $\mathbb{Z}$, with density $\lambda_r$ of particles of class $r$.*

There is an analogous result for the case of $\mathbb{Z}_N$. First we construct operators $F_m^{(N)} : \mathcal{U}_m^{(N)} \times \mathcal{U}_1^{(N)} \mapsto \mathcal{U}_{m+1}^{(N)}$ analogous to the queueing operators $F_m$ above.

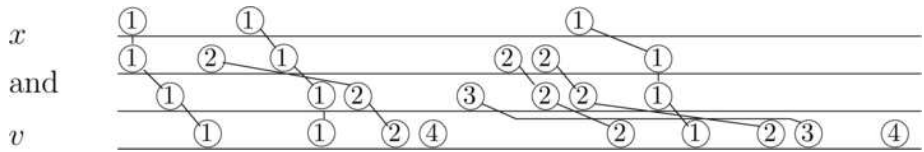

FIG. 3. *Multi-line representation for 3 types of arrivals. Each arrival has been linked to the corresponding service.*

Suppose we have an "arrival process" (with $m$ types of customer and holes) $a \in \mathcal{U}_m^{(N)}$, and a "service process" $s \in \mathcal{U}_1^{(N)}$. We will derive a "departure process" (with $m+1$ types and holes) $d \in \mathcal{U}_{m+1}^{(N)}$ using an extension of the "collapse" described above.

Let $\mathcal{A}_r = \{i \in \mathbb{Z}_N : a(i) = r\}$ for $r \le m$, and let $\mathcal{S} = \{i \in \mathbb{Z}_N : s(i) = 1\}$.

Assume that $|\mathcal{A}_1 \cup \cdots \cup \mathcal{A}_m| \le |\mathcal{S}|$. We will construct the sets

$$\mathcal{D}_r = \{i : d(i) = r\}, \qquad r \le m+1.$$

First $\mathcal{D}_1$ is constructed gradually as follows. $\mathcal{D}_1$ starts empty; take each site $i \in \mathcal{A}_1$ one by one in some arbitrary order, and add to the set $\mathcal{D}_1$ the first element of the list $\{i, i+1, i+2, \ldots\}$ which is in the set $\mathcal{S}$ and which has not yet been added to $\mathcal{D}_1$. (As before, the order in which the sites in $\mathcal{A}_1$ are considered makes no difference to the set $\mathcal{D}_1$ finally obtained.)

Now suppose that $1 < r \le m$ and that the sets $\mathcal{D}_s$ have been constructed for $s < r$. The set $\mathcal{D}_r$ is constructed using the same collapse process, based on the set $\mathcal{A}_r$ and the set $\mathcal{S} \setminus (\mathcal{D}_1 \cup \cdots \cup \mathcal{D}_{r-1})$. Start the set $\mathcal{D}_r$ empty. Take each site $i \in \mathcal{A}_r$ one by one, and add to $\mathcal{D}_r$ the first $j$ in the list $\{i, i+1, i+2, \ldots\}$ such that $j \in \mathcal{S} \setminus (\mathcal{D}_1 \cup \cdots \cup \mathcal{D}_{r-1})$ and such that $j$ has not yet been added to the set $\mathcal{D}_r$.

In this way the sets $\mathcal{D}_1, \ldots, \mathcal{D}_m$ are constructed. Finally let $\mathcal{D}_{m+1} = \mathcal{S} \setminus (\mathcal{D}_1 \cup \cdots \cup \mathcal{D}_m)$.

Note that $|\mathcal{D}_r| = |\mathcal{A}_r|$ for each $r \le m$, and $|\mathcal{D}_{m+1}| = |\mathcal{S}| - |\mathcal{A}_1| - \cdots - |\mathcal{A}_m|$.

We set $d(i) = r$ if $i \in \mathcal{D}_r$, for $r = 1, \ldots, m+1$, and otherwise $d(i) = \infty$.

The operator $F_m^{(N)}$ is then defined by $F_m^{(N)}(a, s) = d$.

Fix $n > 1$. Let $\mathcal{X}^{(N)}$ be the state-space $\{1, \infty\}^{\mathbb{Z}_N \times \{1, \ldots, n\}}$ of configurations of particles with $n$ lines. For $x = (x_m(i), i \in \mathbb{Z}_N, m \in \{1, 2, \ldots, n\}) \in \mathcal{X}^{(N)}$, we say that $x$ has a particle at site $i$ on line $m$ if $x_m(i) = 1$, and that $x$ has a hole at site $i$ on line $m$ if $x_m(i) = \infty$.

Consider a richer state space $\mathcal{V}^{(N)}$ consisting of configurations $v = (v_m(i), i \in \mathbb{Z}_N, m \in \{1, 2, \ldots, n\})$ with $v_m(i) \in \{1, 2, \ldots, m\} \cup \{\infty\}$. If $v_m(i) = \infty$ we say that $v$ has a hole at site $i$ on line $m$, while if $v_m(i) = r < \infty$ then $v$ has a particle of class $r$ at site $i$ on line $m$. Thus $v_m$, the $m$th line of $v$, is a configuration in $\mathcal{U}_m^{(N)}$; all the particles on line $m$ have classes in $\{1, 2, \ldots, m\}$.

Given a configuration $x \in \mathcal{X}^{(N)}$, we will associate a configuration $v = V^{(N)}x \in \mathcal{V}^{(N)}$ to it. We will have $v_m(i) = \infty$ iff $x_m(i) = \infty$; that is, the particles in $v$ and $x$ are in exactly the same places, and deriving $v$ from $x$ corresponds to assigning classes to the particles of $x$.

First set $v_1(\cdot) = x_1(\cdot)$. (All the particles on the top line have class 1.) Now recursively, define $v_{m+1}(\cdot) = F_m^{(N)}(v_m(\cdot), x_{m+1}(\cdot))$ for $m = 1, 2, \ldots, n-1$.

THEOREM 2.2. *Let $p_r \in \mathbb{N}$ for $1 \le r \le n$, with $p_1 + \cdots + p_n \le N$. Let $q_m = p_1 + \cdots + p_m$, for $1 \le m \le n$. Define $\mathcal{X}_{\mathbf{q}}^{(N)}$ to be the subset of $\mathcal{X}^{(N)}$* *consisting of those configurations which have exactly $q_m$ particles on line $m$ for $m = 1, 2, \ldots, n$. Suppose that $x$ has the uniform distribution on $\mathcal{X}_{\mathbf{q}}^{(N)}$. The resulting distribution of the configuration $v_m$ (the mth line of the configuration $v$) is the stationary distribution of the m-type TASEP on $\mathbb{Z}_N$, with $p_r$ particles of class $r$.*

The proofs of Theorems 2.1 and 2.2 are analogous and are done together in the next two sections.

The proof consists of two steps. First we define a Markov process on the state-space $\mathcal{X}$ (or $\mathcal{X}^{(N)}$) (the "multiline process") and show that the product distribution (or the uniform distribution) defined above is a stationary distribution for the process.

Then we show that if $x$ evolves according to the multiline process, then $v = Vx$ (or $V^{(N)}x$) evolves in such a way that its bottom line $(Vx)_n$ realizes a TASEP with $n$ types of particle.

Thus the image of the bottom line of $v$ under the product distribution (or uniform distribution) for $x$ is indeed stationary for the $n$-type TASEP.

2.3. *Further details of multi-class queue construction.* In this section we give a more formal construction of the queue-length process for the multiclass queues described above, and give recurrences analogous to (2.4) for the multiclass case.

First the case of $\mathbb{Z}$. For $k \leq m$, let $A_{[i,j]}^{\leq k}$ be the number of customers arriving in the interval $[i, j]$ whose class is less than or equal to $k$; that is, $A_{[i,j]}^{\leq k} = \#\{r \in [i, j] : a(r) \leq k\}$. Let $S_{[i,j]} = \#\{r \in [i, j] : s(r) = 1\}$ be the number of services occurring in $[i, j]$. We now define

$$(2.7) \qquad Q_j^{\leq k} = \sup_{i \leq j}[A_{[i,j]}^{\leq k} - S_{[i,j]}]_+.$$

The quantity $Q_j^{\leq k}$ represents the number of customers in the queue whose class is less than or equal to $k$, just after time $j$. The definition (2.7) corresponds to (2.3), applied to the queue in which we consider all customers of class less than or equal to $k$ as equivalent, and ignore all other customers. Then, analogously to (2.4),

$$(2.8) \qquad Q_j^{\leq k} = [Q_{j-1}^{\leq k} + I(a(j) \leq k) - I(s(j) = 1)].$$

The interpretation is as follows: at each time step, the number of customers in the queue of class less than or equal to $k$ increases by 1 if such a customer arrives and there is no service, decreases by 1 if no such customer arrives and there is a service (unless it was already zero), and otherwise stays the same.

The departure process $d = F_m(a, s)$ from the queue is given by:

$$
\begin{aligned}
\text{for } k \leq m, \qquad & d(j) \leq k, && \text{if } s(j) = 1 \text{ and either } Q_{j-1}^{\leq k} > 0 \text{ or } a(j) \leq k; \\
(2.9) \qquad & d(j) = m+1, && \text{if } s(j) = 1, Q_{j-1}^{\leq m} = 0 \text{ and } a(j) = \infty; \\
& d(j) = \infty, && \text{if } s(j) = \infty.
\end{aligned}
$$

For $j \in \mathbb{Z}$, write $\overline{Q}_j$ for the vector $(Q_j^{\leq 1}, Q_j^{\leq 2}, \ldots, Q_j^{\leq m})$. Note that since $a(i), s(i), i \in \mathbb{Z}$ are all independent, and the processes $a$ and $s$ are stationary, we have from (2.7) and (2.8) that the queue-length process $\overline{Q}_j, j \in \mathbb{Z}$ is a stationary Markov chain. Each component $Q_j^{\leq k}$ considered in isolation corresponds to a simple queue of the sort described in Section 2.1, with arrival rate $\lambda_1 + \lambda_2 + \cdots + \lambda_k$ and service rate $\mu = \lambda_1 + \cdots + \lambda_{m+1}$.

In the case of $\mathbb{Z}_N$ one defines $A_{[i,j]}$ and $S_{[i,j]}$ in the same way and then

$$
(2.10) \qquad Q_j^{\leq k} = \sup_i [A_{[i,j]}^{\leq k} - S_{[i,j]}]_+;
$$

one again obtains (2.8) and can define $d = F_m^{(N)}(a, s)$ by (2.9). For $j \in \mathbb{Z}_N$, we again write $\overline{Q}_j$ for the vector $(Q_j^{\leq 1}, Q_j^{\leq 2}, \ldots, Q_j^{\leq m})$. It can still be helpful to think of $Q_j^{\leq k}$ as the "number of customers of type less than or equal to $k$ present in the queue just after time $j$" in order to interpret the recurrence (2.8) as the operation of a (priority) queueing server; now this interpretation only holds "locally," since the set of "times" is $\mathbb{Z}_N$ rather than $\mathbb{Z}$.

See Figures 1, 2 and 3 for illustrations of the multiclass and multiline constructions.

## 3. Multiline process.
The multiline process with $n$ lines has state space $\mathcal{X}$ (or $\mathcal{X}^{(N)}) = \{1, \infty\}^{\mathbb{Z} \times \{1, \ldots, n\}}$ (or $\{1, \infty\}^{\mathbb{Z}_N \times \{1, \ldots, n\}}$, resp.). For $m \in \{1, \ldots, n\}$ and $i \in \mathbb{Z}$ or $\mathbb{Z}_N$, we say that $x_m(i)$ is the value at site $i$ on line $m$. The value 1 indicates a particle and the value $\infty$ indicates a hole.

In each of the $n$ lines, the transitions that occur are those of the TASEP (with one type of particle): when a bell rings at site $i$ on line $m$, the configuration $x_m(\cdot)$ moves to $x_m^{(i-1,i)}(\cdot)$. Bells on the different lines all ring together (but not necessarily at the same site).

We have a process of bells at rate 1 at each site $i$ of the bottom line (line $n$). When a bell rings on line $n$ it generates a bell on the line above (line $n-1$), which in turns generates one on the line above that (line $n-2$) and so on up to 1. The rule is as follows. Suppose a bell rings at site $i$ on line $m$. Then if $x_m(i) = 1$ (before any swap), the bell on line $m-1$ also rings at site $i$, while if $x_m(i) = \infty$, then the bell on line $m-1$ rings at site $i+1$.

Formally, we can define as follows. For $x \in \mathcal{X}$ (or $\mathcal{X}^{(N)}$) and $i \in \mathbb{Z}$ (or $\mathbb{Z}_N$), define the sequence $b_m = b_m(x, i), m = 0, \ldots, n$ by $b_n = i$ and, for $m =$

$n-1, \ldots, 0$,

(3.1)
$$b_m = \begin{cases} b_{m+1}, & \text{if } x_{m+1}(b_{m+1}) = 1, \\ b_{m+1} + 1, & \text{if } x_{m+1}(b_{m+1}) = \infty. \end{cases}$$

Thus $b_m$ is the location of the bell which rings on line $m$, for $m = 1, \ldots, n$ (and we have also defined the value $b_0$ for later use). So when the bell rings at site $i$ on line $n$, one jumps from the state $x$ to the state $y = Y(x, i)$ defined by

(3.2)
$$y_m(j) = x_m(j) \qquad \text{for } j \notin (b_m - 1, b_m),$$
$$y_m(b_m - 1) = \min\{x_m(b_m - 1), x_m(b_m)\},$$
$$y_m(b_m) = \max\{x_m(b_m - 1), x_m(b_m)\},$$

where $b_m = b_m(x, i)$ throughout.

The multiline process is illustrated in Figure 4.

THEOREM 3.1. (i) *Let* $q_m \in \{0, 1, \ldots, N\}$ *for* $m \in \{1, 2, \ldots, n\}$. *As before, let* $\mathcal{X}_{\mathbf{q}}^{(N)}$ *be the subset of* $\mathcal{X}^{(N)}$ *consisting of those configurations which have exactly* $q_m$ *particles on line* $m$, *for* $m = 1, 2, \ldots, n$. *Then the uniform distribution on* $\mathcal{X}_{\mathbf{q}}^{(N)}$ *is a stationary distribution for the multiline process on* $\mathbb{Z}^{(N)}$.

(ii) *Let* $\rho_1, \rho_2, \ldots, \rho_n \in [0, 1]$. *Let* $\nu$ *be the product distribution on* $\mathcal{X}$ *under which all the* $\{x_m(i)\}$ *are independent, and* $\mathbb{P}(x_m(i) = 1) = \rho_m$. *Then* $\nu$ *is a stationary distribution for the multiline process on* $\mathbb{Z}$.

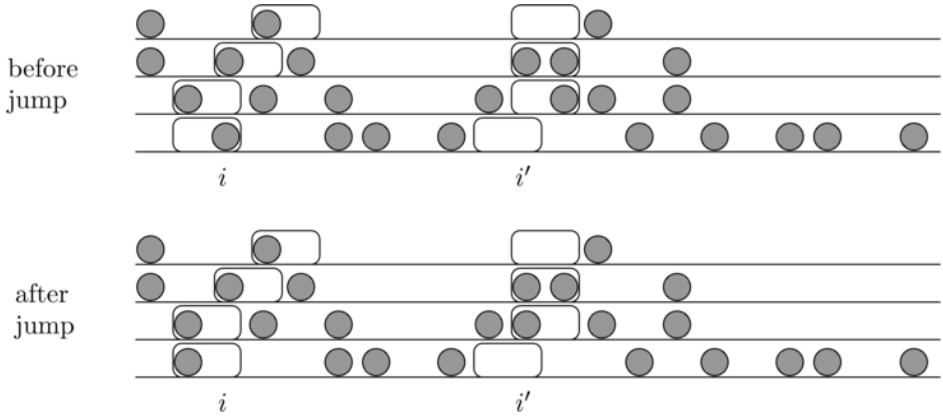

FIG. 4. *Multiline process. In this picture we illustrate the effect of two possible bells, at $i$ and $i'$. For $i$, the bells ring at $b_4 = b_3 = i$, $b_2 = i + 1$ and $b_1 = i + 2$, and $b_0 = i + 3$. The only particle which jumps is on the bottom line, due to the bell at $i$. For $i'$ the bells ring at $b_4' = i'$ and $b_3' = b_2' = b_1' = i' + 1$, and $b_0' = i' + 2$. A particle in the third line jumps due to this bell.*

The following proposition (which is not needed for the main proof) shows that, in equilibrium, the multiline process provides a coupling of $n$ stationary TASEP processes (which, by the previous theorem, have independent marginals at any fixed time). It will follow very easily from our proof of Theorem 3.1; the argument is given at the end of this section.

PROPOSITION 3.2. *Consider the multiline process in one of the stationary distributions described in Theorem* 3.1. *The marginal law of each line of the process is that of a TASEP.*

(In fact, it is possible to go further and show that, in equilibrium, the process of bells on any line of the multiline process consists of a collection of independent Poisson processes. See [6] for details of this extension in the case of $\mathbb{Z}$.)

The rest of this section is devoted to the proof of Theorem 3.1. We will do the proof for (i), the process on $\mathbb{Z}_N$. Part (ii) then follows by letting $N \to \infty$ and choosing $q_1^{(N)}, \ldots, q_n^{(N)}$ such that $q_m^{(N)}/N \to \rho_m$. In that case: (a) the uniform measure of part (i) converges to the product measure of part (ii); (b) one can construct jointly the processes on $\mathbb{Z}_{2N+1}$ and on $\mathbb{Z}$ in such a way that the bells at sites $N+1, \ldots, 2N, 0, 1, \ldots, N$ of the process on $\mathbb{Z}_{2N+1}$ are the same as the bells at sites $-N, \ldots, -1, 0, \ldots, N$ of the process on $\mathbb{Z}$. Combining (a) and (b), one can obtain a coupling of the sequence of finite processes (with initial conditions according to the uniform measure) with the infinite process (with initial condition according to the product measure) with the following property: for any fixed time $t$ and on any fixed finite window $[-k, k]$, the evolution of the processes on $\mathbb{Z}_{2N+1}$ and on $\mathbb{Z}$ over the time interval $[0, t]$ are identical at all sites in $[-k, k]$ for all large enough $N$, with probability 1. In particular, the expected values of any cylinder function of the state at time 0 and at time $t$ for the finite processes converge to the corresponding expected values for the infinite process. Hence since the uniform measures are stationary for the finite processes, one has also that the product measure obtained in the limit is stationary for the infinite process.

So from now on we concentrate on the case of $\mathbb{Z}_N$. For convenience we will write $\tilde{\mathcal{X}} = \mathcal{X}_{\mathbf{q}}^{(N)}$.

We can describe the process as follows. From any state, the process chooses, at rate $N$, a random $i$ uniformly from $\{1, 2, \ldots, N\}$ and jumps from its current state $x$ to the state $Y(x, i)$ defined at (3.2). [Note that we may have $Y(x, i) = x$ so that the process actually stays in its current state.]

Consider the function $T$ from $\tilde{\mathcal{X}} \times \{1, \ldots, N\}$ to itself defined by

(3.3)    $$T(x, i) = (Y(x, i), b_0(x, i)),$$

where $b_0$ is defined at (3.1).

We will show that $T(\cdot, \cdot)$ is a bijective function. From this it follows that the uniform distribution is stationary. To see this, one simply calculates the overall rate of jumps into and out of each state $y$. For convenience we include jumps from $y$ to itself which occur under the description above. Let $\mu(x) = 1/|\tilde{\mathcal{X}}|$ be the probability of a state $x$. The rate of jump out of $y$ is simply $N\mu(y) = N/|\tilde{\mathcal{X}}|$. For the rate of jump into $y$, we have

$$\sum_i \sum_x \mu(x) I(y = Y(x,i)) = \frac{1}{|\tilde{\mathcal{X}}|} \sum_i \sum_x \sum_j I(T(x,i) = (y,j))$$

$$= \frac{1}{|\tilde{\mathcal{X}}|} \sum_j \sum_i \sum_x I((x,i) = T^{-1}(y,j))$$

$$= \frac{1}{|\tilde{\mathcal{X}}|} \sum_j 1$$

$$= N/|\tilde{\mathcal{X}}|.$$

So the rest of the proof is to show that $T$ is indeed a bijective function. We will construct its inverse. In fact, this is essentially equivalent to constructing the time-reversal (in equilibrium) of the process. It turns out that this process can be seen as the image on the original process, reflected between top and bottom and left and right.

In the reverse process, the possible transitions on each line are transitions of a TASEP in which particles move to the right; that is, $(1, \infty)$ may swap to $(\infty, 1)$. A bell at position $j$ on line $m$ will try to swap positions $j-1$ and $j$ on line $m$ in this way. Bells now ring at each site on the top line (line 1) at rate 1, and each bell on line $m$ generates a bell on line $m+1$, as follows. If in state $y$ a bell rings at position $j$ on line $m$, then the bell on line $m+1$ rings at position $j$ if $y_m(j-1) = 1$ (before any swap) and at position $j-1$ if $y_m(j-1) = \infty$.

Formally, for $y \in \tilde{\mathcal{X}}$ and $j \in \mathbb{Z}_N$, define the sequence $c_m = c_m(y,j), m = 1, \ldots, n+1$ by $c_1 = j$ and, for $m = 1, 2, \ldots, n$,

$$(3.4) \qquad c_{m+1} = \begin{cases} c_m, & \text{if } y_m(c_m - 1) = 1, \\ c_m - 1, & \text{if } y_m(c_m - 1) = \infty. \end{cases}$$

Thus when the bell rings at site $j$ on line 1, the bells ring also at sites $c_m(y,j)$ on lines $m = 2, \ldots, n$, and the process jumps from the state $y$ to the state $x = X(y,j)$ defined by

$$(3.5) \qquad x_m(i) = y_m(i) \text{ for } i \notin (c_m - 1, c_m),$$

$$(3.6) \qquad x_m(c_m - 1) = \max\{y_m(c_m - 1), y_m(c_m)\},$$

$$(3.7) \qquad x_m(c_m) = \min\{y_m(c_m - 1), y_m(c_m)\},$$

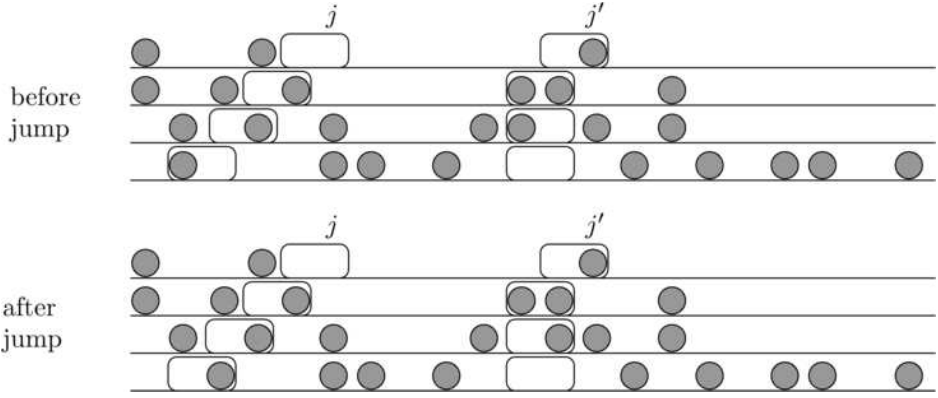

Fig. 5.    *Reverse multiline process. We illustrate the two reverse jumps which correspond to the forward jumps in Figure 4. The corresponding reversed bells ring at sites $j$ and $j'$. For $j$, the bells ring at $c_1 = j$, $c_2 = j-1$, $c_3 = j-2$ and $c_4 = j-3$; also $c_5 = j-3$. Again, only one particle, on the bottom line, jumps due to the bell at $j$. For $j'$ the bells ring at $c_1' = j'$ and $c_2' = c_3' = c_4' = j'-1$, with $c_5' = j'-2$. Again, a particle in the third line jumps due to this bell.*

where $c_m = c_m(y, j)$ throughout. We have also defined $c_{n+1}$.

This reverse process is illustrated in Figure 5.

Now define the function $T^*$ from $\tilde{\mathcal{X}} \times \{1, \ldots, N\}$ to itself by

$$(3.8) \qquad T^*(y, j) = (X(y, j), c_{n+1}(y, j)).$$

This is the inverse of $T$ that we need:

PROPOSITION 3.3.    *If $x \in \tilde{\mathcal{X}}$ and $i \in \{1, \ldots, N\}$, then $T^*(T(x, i)) = (x, i)$.*

PROOF.    Let $(y, j) = T(x, i)$, so that $y = Y(x, i)$ and $j = b_0(x, i)$. Write $b_m = b_m(x, i)$ and $c_m = c_m(y, j)$.

Then the result we need is that $X(y, j) = x$ and $c_{n+1} = i$.

The configurations $x$ and $y$ are identical except in the following way: for each $m$, if $(x_m(b_m - 1), x_m(b_m)) = (\infty, 1)$ then these values are swapped to $(1, \infty)$ in $y$.

Similarly, the configurations $y$ and $X(y, j)$ are identical except as follows: for each $m$, if $(y_m(c_m - 1), y_m(c_m)) = (1, \infty)$ then these values are swapped to $(\infty, 1)$ in $X(y, j)$.

So to show that $X(y, j) = x$, it is enough to show that, for each $m$,

$$(3.9) \ (x_m(b_m - 1), x_m(b_m)) = (\infty, 1) \ \Leftrightarrow \ (y_m(c_m - 1), y_m(c_m)) = (1, \infty),$$

and that in this case $b_m = c_m$.

We first note a few simple implications for later reference. For $m = 1, \ldots, n$:

$$(3.10) \qquad x_m(b_m) = 1 \Rightarrow y_m(b_m - 1) = 1,$$

$$(3.11) \qquad x_m(b_m) = \infty \Rightarrow y_m(b_m) = \infty,$$

$$(3.12) \quad (x_m(b_m - 1), x_m(b_m)) = (\infty, 1) \Rightarrow (y_m(b_m - 1), y_m(b_m)) = (1, \infty),$$

$$(3.13) \quad x_m(b_m) = 1 \text{ and } y_m(b_m) = \infty \Rightarrow x_m(b_m - 1) = \infty.$$

These all follow immediately from the relationship between $x$ and $y$ observed above.

Now, let $\mathcal{C}(m)$ be the property that

$$(3.14) \qquad c_m = b_{m-1} = \begin{cases} b_m, & \text{if } x_m(b_m) = 1, \\ b_m + 1, & \text{if } x_m(b_m) = \infty. \end{cases}$$

(The second equality is the definition of $b_{m-1}$.) We will show by induction that $\mathcal{C}(m)$ holds for $m = 1, \ldots, n, n+1$.

The property $\mathcal{C}(1)$ is true by definition, since $c_1 = c_1(y, b_0) = b_0$.

Suppose $r \in \{1, \ldots, n\}$ and $\mathcal{C}(r)$ holds, so that

$$(3.15) \qquad c_r = b_{r-1} = \begin{cases} b_r, & \text{if } x_r(b_r) = 1, \\ b_r + 1, & \text{if } x_r(b_r) = \infty. \end{cases}$$

There are two cases to check:

*Case* 1. $x_r(b_r) = 1$. Then $c_r = b_r$ [from (3.15)], and also $y_r(b_r - 1) = 1$, from (3.10). So $y_r(c_r - 1) = 1$, and $c_{r+1} = c_r = b_r$, using (3.4). Hence $\mathcal{C}(r+1)$ holds as desired.

*Case* 2. $x_r(b_r) = \infty$. Then $c_r = b_r + 1$ [from (3.15)], and $y_r(b_r) = \infty$ [from (3.11)]. So $y_r(c_r - 1) = y_r(b_r) = \infty$, and (3.4) gives $c_{r+1} = c_r - 1 = b_r$, giving $\mathcal{C}(r+1)$ as desired.

Thus by induction $\mathcal{C}(r)$ holds for all $r = 1, 2, \ldots, n+1$.

In particular putting $r = n + 1$ we have that $c_{n+1} = b_n = i$. So it remains to show that the two sides of (3.9) are equivalent, and that each implies that $c_m = b_m$.

First suppose that the left-hand side of (3.9) holds, that is,

$$(x_m(b_m - 1), x_m(b_m)) = (\infty, 1).$$

We have $c_m = b_m$ by the property (3.14), and then, using (3.12),

$$(y_m(c_m - 1), y_m(c_m)) = (y_m(b_m - 1), y_m(b_m))$$
$$= (1, \infty),$$

which is the right-hand side of (3.9). In the other direction, suppose that the right-hand side of (3.9) holds, that is,

$$(y_m(c_m - 1), y_m(c_m)) = (1, \infty).$$

Suppose that we also had $x_m(b_m) = \infty$. Then we would have $y_m(b_m) = \infty$ from (3.11). Also from (3.14) we would have $c_m = b_m + 1$, so that $y_m(b_m) = y_m(c_m - 1) = 1$, which is a contradiction.

So we must in fact have $x_m(b_m) = 1$. Then from (3.14), $c_m = b_m$ as desired. We have $y_m(b_m) = y_m(c_m) = \infty$, so from (3.13) we have $x_m(b_m - 1) = \infty$. Thus the left-hand side of (3.9) indeed holds as desired. This completes the proof of Proposition 3.3 and hence of Theorem 3.1. □

PROOF OF PROPOSITION 3.2. By construction, line 1 of the reverse process realizes a (reverse) TASEP. Hence line 1 of the forward process realizes a forward TASEP. Now let $1 < m \leq n$ and consider the process consisting just of lines $m, m+1, \ldots, n$. This is a multiline process with $n - m + 1$ lines. For exactly the same reason, the top line of this process realizes a TASEP. But this is exactly line $m$ of the original process. So (in equilibrium) every line realizes a TASEP, as desired. □

## 4. Realization of the multitype TASEP by the multiline process.

THEOREM 4.1 (for both $\mathbb{Z}$ and $\mathbb{Z}_N$). *If $x$ evolves according to the multiline process (with $n$ lines) then $v = Vx$ (or $V^{(N)}x$) evolves in such a way that its $n$th line $v_n$ follows a TASEP with $n$ types.*

Since the bells at each site of level $n$ of the multiline process occur at rate 1 independently, it will be enough to show the following:

PROPOSITION 4.2. *Let $m \in \{1, 2, \ldots, n\}$. When a bell rings at site $i$ on line $m$ of the multiline process, the effect on $v_m$ is also of a bell at site $i$, that is, $v_m$ jumps to $v_m^{(i-1,i)}$.*

This is clearly true for $m = 1$, since $v_1 = x_1$.

For $m \geq 2$, recall that a bell at site $i$ on level $m$ creates a bell on level $m - 1$ according to the following rule: if there was a particle at site $i$ on level $m$, then the bell on level $m - 1$ is also at site $i$; if not, then the bell on level $m - 1$ is at site $i + 1$.

Using this, we can prove the proposition by induction. Under the assumption that the bell on level $m - 1$ has the effect of a TASEP bell, we show that the swaps on level $m - 1$ and $m$ also have the effect of a TASEP bell on level $m$. In terms of the arrival, service and departure processes, the induction step is the following lemma; to perform the induction, we apply it with $a = v_{m-1}$, $s = x_m$ and $d = v_m$.

LEMMA 4.3.   (i) *Let $m \geq 2$. Let $a$ be an arrival process with $m-1$ types, and let $s$ be a service process. Let $d = F_m(a, s)$ be the corresponding departure process with $m$ types. Define*

$$b = \begin{cases} i, & \text{if } s(i) = 1, \\ i+1, & \text{if } s(i) = \infty. \end{cases}$$

*Write $\tilde{a} = a^{(b-1,b)}$, $\tilde{s} = s^{(i-1,i)}$ and $\tilde{d} = F_m(\tilde{a}, \tilde{s})$. Then $\tilde{d} = d^{(i-1,i)}$.*

(ii) *The analogous statement for the process on $\mathbb{Z}_N$ rather than $\mathbb{Z}$; the operator $F_m$ is replaced by $F_m^{(N)}$.*

PROOF.   We will word the proof for the case of $\mathbb{Z}$; but an identical argument applies to the case of $\mathbb{Z}_N$, using the framework explained at the end of Section 2.3. We will check three cases:

*Case* 1. $s(i) = \infty$. There is no service at $i$, and the bell at $i$ leaves the process $s$ unchanged: $\tilde{s} = s$. Also $d(i) = \infty$, so $d^{(i-1,i)} = d$. We have $b = i + 1$.

The order of events in the arrival process at $i$ and $i + 1$ does not matter; since there is no service at $i$, any arrivals at $i$ or $i + 1$ have no opportunity to depart until $i + 1$ or later. So $\tilde{d} = F_m(\tilde{a}, \tilde{s}) = F_m(\tilde{a}, s) = F_m(a, s) = d = d^{(i-1,i)}$, as required.

*Case* 2. $s(i-1) = \infty$ and $s(i) = 1$. Thus $\tilde{s}(i-1) = s(i)$ and $\tilde{s}(i) = \infty$. Similarly $d^{(i-1,i)}(i-1) = d(i)$ and $d^{(i-1,i)}(i) = d(i-1) = \infty$. We have $b = i$. Let $\mathcal{C}$ be the collection of all the customers who were either present in the queue just after time $i - 2$ or who arrived at time $i - 1$ or $i$. The value of $d(i)$ is the class of the highest priority customer in the collection $\mathcal{C}$ (or is equal to $m + 1$ if $\mathcal{C}$ is empty).

After the swaps, the higher priority out of the particles arriving at $(i-1)$ and $i$ (if any) has moved to place $(i-1)$ in the arrival process. So the particle which previously departed at time $i$ (if any) is now already present at time $i - 1$. Since $\tilde{s}(i-1) = 1$, this particle will depart at time $i - 1$. That is, $\tilde{d}(i-1) = d(i)$, and since $\tilde{s}(i) = \infty$, we have also $\tilde{d}(i) = \infty = d(i-1)$. Overall, the collections of particles arriving and departing in $[i-1, i]$ are unchanged and so $\overline{Q}_i$, and hence all later values in the queue-length process and the departure process, are unchanged. Thus $\tilde{d} = d^{(i-1,i)}$ as desired.

*Case* 3. $s(i-1) = s(i) = 1$. There are services available at both $i - 1$ and $i$. Consider the collection $\mathcal{C}$ as in case 2 above. The values of $d(i-1)$ and $d(i)$ are, in some order, the highest and second highest priorities out of this collection $\mathcal{C}$. [The two values may be the same. Also, if this collection has fewer than 2 particles, then one or both of $d(i-1)$ and $d(i)$ will be equal to $m + 1$.]

The bell has no effect on the service process: $\tilde{s} = s$. For the arrival process, we have $b = i$, and the swap at $(i-1, i)$ ensures that $\tilde{a}(i-1) \leq \tilde{a}(i)$, so that the highest priority customer out of the collection $\mathcal{C}$ is available to depart at time $i-1$ (since this customer is either already present in the queue after time $i-2$ or arrives at time $i-1$). Thus this customer (which was the higher priority of those departing at times $i-1$ and $i$ in the process $d$) will now depart at $i-1$. The second-highest now departs at $i$. Again $\overline{Q}_i$ is left unchanged. Thus $\tilde{d} = d^{(i-1,i)}$. $\quad\square$

## 5. Common denominators and states of minimal weight for the process on $\mathbb{Z}_N$.

THEOREM 5.1. *Consider the $n$-type process on $\mathbb{Z}_N$, with $p_m$ particles of type $m$, for $m = 1, 2, \ldots, n$, where $p_1 + p_2 + \cdots + p_n \leq N$.*

(i) *In the stationary distribution, the probability of any state is an integer multiple of $M^{-1}$, where*

$$(5.1) \qquad M = \left[ \binom{N}{p_1} \binom{N}{p_1 + p_2} \cdots \binom{N}{p_1 + p_2 + \cdots + p_n} \right].$$

(ii) *A state $u$ has probability $M^{-1}$ (the smallest possible) in the stationary distribution iff the following holds for each $j \in \mathbb{Z}_N$: if $u(j) = \infty$ then $u(j+1) = \infty$ or $u(j+1) = n$, while if $u(j) = m < \infty$ then $u(j+1) \geq m-1$.*

We write $\mathcal{W}_n^{(N)}$ for the set of states satisfying the condition of part (ii). The proof of the first part is immediate:

PROOF OF THEOREM 5.1(i). $M$ is the number of configurations of the multiline process, with $(p_1 + \cdots + p_m)$ particles on line $m$. By Theorem 2.2, the stationary distribution for the $n$-TASEP is given by the image, under the map $x \to (V^{(N)}x)_n$, of the uniform distribution for the multiline process. Thus the probability of any state is a multiple of $M^{-1}$, as required for part (i). $\quad\square$

Before proving the second part, we need a couple of definitions and lemmas. Recall the "queueing operator" $F_m^{(N)} : \mathcal{U}_m^{(N)} \times \mathcal{U}_1^{(N)} \mapsto \mathcal{U}_{m+1}^{(N)}$. The next two results concern solutions of the relation $d = F_m^{(N)}(a, s)$. Here $a \in \mathcal{U}_m^{(N)}$ is an "arrival process" (consisting of $m$ classes of arrival, corresponding to values $a(j) \in \{1, 2, \ldots, m\}$, and holes, corresponding to the value $\infty$); $s \in \mathcal{U}_1^{(N)}$ is a "service process" (consisting of values 1, where service is possible, and values $\infty$, where no service is possible), and $d \in \mathcal{U}_{m+1}^{(N)}$ is a "departure process" (consisting of $m$ classes of departure, corresponding to values in $\{1, 2, \ldots, m\}$, of "unused services," corresponding to the value $m+1$, and holes or "nonservices," corresponding to the value $\infty$).

The first lemma shows in effect that $s$ is uniquely determined by $d$, and that there one can always find at least one appropriate $a$:

LEMMA 5.2.  *Let $d \in \mathcal{U}_{m+1}^{(N)}$.*

(i)  *Define $H(d) \in \mathcal{U}_1^{(N)}$ by*

$$(H(d))(j) = \begin{cases} \infty, & \text{if } d(j) = \infty, \\ 1, & \text{if } d(j) < \infty. \end{cases}$$

*If $F_m^{(N)}(a, s) = d$ for some $a, s$, then $s = H(d)$.*

(ii)  *Define further $G_m(d) \in \mathcal{U}_m^{(N)}$ by*

$$(G_m(d))(j) = \begin{cases} \infty, & \text{if } d(j) = m + 1 \text{ or } d(j) = \infty, \\ d(j), & \text{if } d(j) \le m. \end{cases}$$

*Then $F_m^{(N)}(G_m(d), H(d)) = d$.*

PROOF.  Part (i) follows directly from the definition of $F_m^{(N)}$, since if $d = F_m^{(N)}(a, s)$, then $d(j) = \infty$ iff $s(j) = \infty$. Part (ii) is easy to check. In terms of the collapse process, there is a space available directly below the starting position of every particle; so each particle remains at the same site where it starts. In terms of the queue, there is a service at every time when an arrival occurs; thus each arrival departs immediately and the queue is always empty.  □

The following lemma now shows that $a = G_m(d)$ gives the *unique* solution whenever $d \in \mathcal{W}_{m+1}^{(N)}$, and that the property defined in Theorem 5.1(ii) is inherited by $G_m(d)$:

LEMMA 5.3.  *If $d \in \mathcal{W}_{m+1}^{(N)}$ and $F_m^{(N)}(a, s) = d$, then $a = G_m(d)$ and $s = H(d)$. Furthermore, $G_m(d) \in \mathcal{W}_m^{(N)}$.*

PROOF.  The fact that $s = H(d)$ is already given by the Lemma 5.2(i), so it remains to look at the process $a$.

Divide the configuration $d$ into blocks, each of which consists of a string of holes followed by a string of particles. Suppose that the block occupies the sites $i_0, \ldots, i_1, \ldots, i_2$, with $i_0 \le i_1 < i_2$, that $d(j) = \infty$ for $j \in [i_0, i_1]$, and that $d(j) < \infty$ for $j \in [i_1 + 1, i_2]$. Note that the definition of $\mathcal{W}_{m+1}^{(N)}$ implies that $d(i_1 + 1) = m + 1$.

Since there is an "unused service" at $i_1 + 1$, and no services in $[i_0, i_1]$, there can be no "arrivals" in $[i_0, i_1]$; that is, $a(j) = \infty$ for $j \in [i_0, i_1]$; and the "queue" is empty after time $i_1$. Since there are no more "nonservices" in $[i_1 + 1, i_2]$, every customer arriving during that time departs immediately. Thus for $j \in [i_1 + 1, i_2]$, we have $a(j) = d(j)$ if $d(j) \le m$, and $a(j) = \infty$ if $d(j) = m + 1$ (unused service).

Hence the block $(a(i_0), \ldots, a(i_2))$ is identical to the string $(d(i_0), \ldots, d(i_2))$, except that any value $d(j) = m+1$ is replaced by the value $a(j) = \infty$. Putting the blocks together, one has exactly that $a = G_m(d)$ as asserted.

Finally, recall that for a configuration in $\mathcal{W}_{m+1}^{(N)}$ the set of "allowable pairs" of values for sites $j$ and $j+1$ are $(\infty, \infty)$, $(\infty, m+1)$ or $(r, k)$ for some $k \geq r-1$. It's easy to check that replacing all the values $m+1$ by $\infty$ transforms any such allowable pair for $\mathcal{W}_{m+1}^{(N)}$ into an allowable pair for $\mathcal{W}_m^{(N)}$. Hence $G_m(d) \in \mathcal{W}_m^{(N)}$ as required.  $\square$

PROOF OF THEOREM 5.1(ii).  A state $u$ has probability exactly $M^{-1}$ iff there is a unique $x$ such that $v_n = u$, where $v = V^{(N)}x$. From the definition of $V^{(N)}$ this means that $v_{m+1} = F_m^{(N)}(v_m, x_{m+1})$ for each $m = 1, 2, \ldots, n-1$, and $v_1 = x_1$.

First, we show how to construct one such $x$ directly. Starting from $u = v_n$, we can recursively set $v_m = G_m(v_{m+1})$ for $m = n-1, \ldots, 1$, and $x_m = H(v_m)$ for each $m$. Lemma 5.2(ii) then gives that indeed $v = V^{(N)}x$.

The configuration $v_n$ is a "minimal state" if this is the *only* possible $x$. It now follows, by applying Lemma 5.3 repeatedly that this is indeed the case whenever $v_n \in \mathcal{W}_n^{(N)}$, since then we can construct $v_m$ uniquely from $v_{m+1}$ for each $m$, and $x_m$ uniquely from $v_m$.

It now remains to show the converse: if $u \notin \mathcal{W}_n^{(N)}$, then there are at least two $x$ such that $(V^{(N)}x)_n = u$. It's enough to show the following: if we are given $v_{m+1} \notin \mathcal{W}_{m+1}^{(N)}$ and are looking for $v_m$ such that $F_m^{(N)}(v_m, x_{m+1}) = v_{m+1}$, then *either* (a) there are two choices for $v_m$ *or* (b) $m > 1$ and the unique choice for $v_m$ gives $v_m \notin \mathcal{W}_m^{(N)}$. This way, starting from $u \notin \mathcal{W}_n^{(N)}$ we can construct all of $v$ (and hence $x$) recursively starting from $v_n = u$, and be guaranteed to have a nonunique choice at some stage.

So, choose $v_{m+1} \notin \mathcal{W}_{m+1}^{(N)}$. Choose some $j$ such that the pair $(v_{m+1}(j-1), v_{m+1}(j))$ violates the conditions defining $\mathcal{W}_{m+1}^{(N)}$.

*Case* 1. $v_{m+1}(j-1) = \infty$ and $v_{m+1}(j) \leq m$. No service occurs at $j-1$, and a customer departs at $j$. Then $v_m = G_m(v_{m+1})$ satisfies $v_m(j-1) = \infty$ and $v_m(j) \leq m$ also. But swapping the values $v_m(j-1)$ and $v_m(j)$ in the arrival process will have no effect on the departure process (since there is no service at $j-1$). Thus the choice of arrival process $v_m$ is nonunique, and (a) holds.

*Case* 2. $v_{m+1}(j) = k \leq m+1$ and $v_{m+1}(j+1) = r \leq k-2$. Certainly $m > 1$. Let $v_m = G_m(v_{m+1})$. Then $v_m(j) \in \{k, \infty\}$, and $v_m(j+1) \leq k-2 \leq m-1$. So $v_m \notin \mathcal{W}_m^{(N)}$, and (b) holds.  $\square$

**6. Independence and comparison properties for the process on $\mathbb{Z}$.** Many important properties of the 2-TASEP on $\mathbb{Z}$ concern the process "as seen from a second-class" particle.

One central result, proved in [4] and [12] may be stated as follows. Let $u(j), j \in \mathbb{Z}$ be a sample from the stationary distribution (with some densities $\rho_1$ and $\rho_2$ of first- and second-class particles). Then conditioned on $u(0) = 2$, the configurations $(u(j), j < 0)$ and $(u(j), j > 0)$ are independent. It is sometimes said that the configuration "factorizes" around the position of a second-class particle. Put another way, the state 2 is a *renewal state* for the process $u(j)$.

As observed by Angel [1], there is a simple proof of this fact via the collapse process, or, equivalently, via the queueing process. Recall that the stationary measure of the 2-TASEP is realized by the departure and service process of a single queue. As explained in Section 2, the queueing system is Markovian; in particular, given the current queue-length, the future and past evolutions of the system are independent (as before, a "time" in the queue corresponds to a site of $\mathbb{Z}$ for the particle system—the "future" and "past" correspond to sites to the right and left respectively). The state $u(0) = 2$ corresponds to an unused service at time 0, and therefore implies that the queue is empty just after time 0; this information "decouples" the past and future as desired.

For the $n$-TASEP with $n \geq 3$, there is no such single renewal state. However, one can find longer "renewal strings." Let $r > 0$ and let $(w(0), w(1), \ldots, w(r))$ be a sequence taking values in $\{1, 2, \ldots, n\}$ satisfying the following conditions:

(i) $w(0) = n$;
(ii) $w(r) = 2$;
(iii) if $n \geq 4$: for each $m \in \{3, \ldots, n-1\}$, there exists $j_m < r$ such that $w(j_m) = m$ and such that $w(j) \leq m$ for all $j \in \{j_m + 1, \ldots, r\}$.

Call such a sequence a renewal string for the $n$-TASEP. For example, for $n = 4$, the string $(4, 1, 2, 3, 1, 2)$ qualifies, with $r = 5$ and $j_3 = 3$. One can show the following:

PROPOSITION 6.1. *Let $w$ be such a renewal string. Conditional on $u(i) = w(i)$ for $0 \leq i \leq r$, the configurations $(u(j), j < 0)$ and $(u(j), j > r)$ are independent.*

PROOF. We do not give the proof in full, but the outline is as follows. Recall that we can write $u = (Vx)_n$, where $x$ is a sample from the stationary distribution of the appropriate multiline process; namely, $x$ is a collection of $n$ Bernoulli processes with appropriate rates. We defined a collection of $n - 1$ queueing processes; the $m$th queue had an arrival process $v_m$ (with $m-1$ types of customer), service process $x_{m+1}$, and departure process $v_{m+1}$ (with $m$ types of customer). Let $\mathcal{C}_j^{(m)}$ denote the collection of customers present in queue $m$ after time $j$. Since the queue-length processes form a Markov chain as described in Section 2, one has that, given the vector of queue-lengths $(\mathcal{C}_j^{(1)}, \ldots, \mathcal{C}_j^{(n-1)})$ after time $j$, the past and future evolution of the process $v$ (and in particular its $n$th line) are independent. Now, the key step is to show that if $(u(0), \ldots, u(r))$ is a renewal string, then $\mathcal{C}_r^{(m)}$ is empty for all $m \in \{1, 2, \ldots, n-1\}$. This can be done by arguments similar to those used in the proof of Theorem 5.1; one constructs $v_m$, $m = n-1, \ldots, 1$ recursively starting from $v_n$, and uses the information that if a particle of type $m$ occurs on line $m$ at site $j$, then the corresponding $(m-1)$st queue is empty immediately after $j$. (These arguments are still valid for the case of $\mathbb{Z}$, since the proofs of Lemma 5.2 and Lemma 5.3 do not use the cyclic structure of $\mathbb{Z}_N$.) From this one obtains that $\mathcal{C}_j^{(m)}$ is empty for all $j_m \leq j \leq r$, where $j_3, \ldots, j_{n-1}$ are defined above and we take $j_2 = r$, $j_n = 0$. $\square$

A related property for the 2-TASEP, stated in [2] and proved in [12] and [4] may be stated as follows: the set $\{j > 0 : u(j) = \infty\}$ is independent of the event $\{u(0) = 2\}$. Similarly, by symmetry, the set $\{j < 0 : u(j) = 1\}$ is independent of the event $\{u(0) = 2\}$. (The relevant symmetry is the following: if one reverses left and right for an $n$-TASEP and also reverses the order of the states $\{1, 2, \ldots, n, \infty\}$, then one again obtains an $n$-TASEP, with the densities appropriately reversed.)

The following version provides a stronger independence property (and also extends to $n > 2$); it is an immediate consequence of the queueing representation.

PROPOSITION 6.2.  *Let $u$ be a configuration distributed according to a translation-invariant stationary distribution for the $n$-TASEP on $\mathbb{Z}$. Then the set $\{j > 0 : u(j) = \infty\}$ is independent of the entire configuration $(u(j) : j \leq 0)$. Similarly, the set $\{j < 0 : u(j) = 1\}$ is independent of the configuration $(u(j), j \geq 0)$.*

PROOF.  Write $u = (Vx)_n$ as above. The configuration $(u(j), j \leq 0)$ depends only on $(x_m(j), j < 0, m \in \{1, \ldots, n\})$, and so is independent of the set $\{j > 0 : u(j) = \infty\} = \{j > 0 : x_n(j) = \infty\}$ (since the state $x$ consists of $n$ independent Bernoulli processes). $\square$

In [2, 4] and [12], the densities of particles and holes around a second class particle are explored. It is shown that for all $j > 0$

$$\mathbb{P}(u(j) = 1) < \mathbb{P}(u(j) = 1 \mid u(0) = 2)$$

and similarly (by symmetry), for all $j < 0$

$$\mathbb{P}(u(j) = \infty) < \mathbb{P}(u(j) = \infty | u(0) = 2).$$

Using the queueing representation, we can extend this to a pathwise stochastic comparison of the stationary measure with the measure conditioned to have a second-class particle at the origin:

PROPOSITION 6.3. *Let $u$ be distributed according to a stationary distribution for the 2-TASEP on $\mathbb{Z}$, with positive density of first- and second-class particles. Let $u'$ be drawn from the same distribution but conditional on the presence of a second-class particle at site $0$. Then there is a coupling of $u$ and $u'$ such that, with probability $1$, one has $u(j) \le u'(j)$ for all $j > 0$ and $u(j) \ge u'(j)$ for all $j < 0$.*

PROOF. Let $\mu < 1$ be the rate of the service process $\mathcal{S}$ and $\lambda < \mu$ the rate of the arrival process $\mathcal{A}$, used to construct the stationary distribution of the 2-TASEP in Section 2. We have

$$Q_j = [Q_{j-1} + I(j \in \mathcal{A}) - I(j \in \mathcal{S})]_+,$$

and if we set

$$u(j) = \begin{cases} 1, & \text{if } j \in \mathcal{S} \text{ and either } Q_{j-1} > 0 \text{ or } j \in \mathcal{A}, \\ 2, & \text{if } j \in \mathcal{S}, Q_{j-1} = 0 \text{ and } j \notin \mathcal{A}, \\ \infty, & \text{if } j \notin \mathcal{S}, \end{cases}$$

then $u$ is distributed according to the stationary distribution of the 2-TASEP with densities $\lambda$ and $\mu - \lambda$ of first- and second-class particles.

Note that $Q_0$ depends only on $\mathcal{A} \cap (-\infty, 0]$ and $\mathcal{S} \cap (-\infty, 0]$; since $\mathcal{A}$ and $\mathcal{S}$ are Bernoulli processes, $Q_0$ is independent of $\mathcal{A} \cap [1, \infty)$ and $\mathcal{S} \cap [1, \infty)$.

Let $Q_0$ be drawn from the stationary queue-length distribution for the queue (which is geometric with parameter $\lambda/\mu$).

Under the condition $u(0) = 2$ we have that the queue-length at time $0$ is $0$. So let $Q'_0 = 0$ and let $Q'_j$ evolve according to

$$Q'_j = [Q'_{j-1} + I(j \in \mathcal{A}) - I(j \in \mathcal{S})]_+.$$

Set also

$$u'(j) = \begin{cases} 1, & \text{if } j \in \mathcal{S} \text{ and either } Q'_{j-1} > 0 \text{ or } j \in \mathcal{A}, \\ 2, & \text{if } j \in \mathcal{S}, Q'_{j-1} = 0 \text{ and } j \notin \mathcal{A}, \\ \infty, & \text{if } j \notin \mathcal{S}. \end{cases}$$

Then $(u'(j), j > 0)$ follows the stationary distribution conditioned on the presence of a second-class particle at $0$, as desired. Comparing the two sets of recurrences, one has $Q_j \ge Q'_j$ for all $j > 0$. The definitions of $u$ and $u'$ then give $u(j) \le u'(j)$ for all $j > 0$ as desired.

The same argument gives the symmetric statement about the processes on $j < 0$. Since the processes for positive and negative $j$ are independent given $Q_0$, the two couplings can in fact be carried out simultaneously as claimed. □

Note that if $T = \inf\{j \geq 0 : Q_j = 0\}$ then $Q_j = Q'_j$ and $u(j) = u'(j)$ for all $j \geq T$. The distribution of $T$ decays exponentially; thus so does the probability that $u(j) \neq u'(j)$. Explicit representations for these probabilities can be found in [2].

We conclude with one more remark concerning the representation of the stationary measure of the 2-TASEP by the $M/M/1$ queueing process. We assigned first-class particles to the sites of departures and second-class particles to the sites of unused services. Using the reversibility of the queue-length process, one can obtain that the joint process of arrivals and unused services is the time-reversal of the joint process of departures and unused services. Thus one could equally assign first-class particles to the sites of arrivals and second-class particles to the sites of unused services, to obtain the stationary measure for a 2-TASEP with jumps to the right rather than to the left. However, it is less clear how to extend this representation naturally to the case $n > 2$.

**Acknowledgments.** We thank the referee for many comments and suggestions which have substantially improved the presentation of the paper. We are grateful for the support of the Fundação de Amparo à Pesquisa do Estado de São Paulo (FAPESP), the Conselho Nacional de Desenvolvimento Científico (CNPq) and the Brazil–France Agreement in Mathematics.

Instituto de Matemática e Estatística
Universidade de São Paulo
Caixa Postal 66281
05311-970 São Paulo
Brazil
E-mail: pablo@ime.usp.br
URL: http://www.ime.usp.br/~pablo

Department of Statistics
University of Oxford
1 South Parks Road
Oxford OX1 3TG
UK
E-mail: martin@stats.ox.ac.uk
URL: http://www.stats.ox.ac.uk/~martin